\documentclass[a4paper,11pt]{article}
\usepackage{amssymb}
\usepackage{amsfonts}
\usepackage{amsmath}
\usepackage{hyperref}
\usepackage{geometry}

\input{tcilatex}
\begin{document}

\title{Stability analysis of an implicit and explicit numerical method for
Volterra integro-differential equations with kernel $K\left( x,y\left(
t\right) ,t\right) $}
\author{J. S. C. Prentice \\
Senior Research Officer\\
Mathsophical Ltd.\\
Johannesburg, South Africa\\
Email: jpmsro@mathsophical.com}
\maketitle

\begin{abstract}
We present implicit and explicit versions of a numerical algorithm for
solving a Volterra integro-differential equation. These algorithms are an
extension of our previous work, and cater for a kernel of general form. We
use an appropriate test equation to study the stability of both algorithms,,
numerically deriving stability regions. The region for the implicit method
appears to be unbounded, while the explicit has a bounded region close to
the origin. We perform a few calculations to demonstrate our results.
\end{abstract}

\section{Introduction}

Recently, we described a numerical method for solving the Volterra
integro-differential equation 
\begin{equation}
y^{\left( n\right) }\left( x\right) =f\left( x,y\right)
+\int\limits_{x_{0}}^{x}Kdt,\text{ \ \ \ }x>x_{0}  \label{problem}
\end{equation}%
with various specific forms for the kernel $K$ \cite{prentice 2}. In this
paper,\ a sequel to \cite{prentice 2}, we consider the more general kernel%
\begin{equation*}
K=K\left( x,y\left( t\right) ,t\right) .
\end{equation*}%
We modify our previous algorithm appropriately, and we consider both
implicit and explicit forms of the resulting algorithm. Our attention will
primarily be focussed on the stability of these algorithms.

\section{Algorithm}

To begin with, we consider the case of $n=1$ in (\ref{problem}). We describe
the more general case in Appendix B. We partition the interval of interest,
denoted $\left[ x_{0},x_{N}\right] ,$ by means of the \textit{equispaced}
nodes%
\begin{equation}
x_{0}<x_{1}<x_{2}<\ldots <x_{N}.  \label{nodes}
\end{equation}%
The spacing between the nodes - the \textit{stepsize} - is denoted $h$.

Using the initial value%
\begin{equation*}
y\left( x_{0}\right) =y_{0},
\end{equation*}%
we compute the solution at $x_{1}$ via%
\begin{align*}
y_{1}& =y_{0}+hf\left( x_{1},y_{1}\right)
+h\int\limits_{x_{0}}^{x_{1}}K\left( x_{1},y\left( t\right) ,t\right) dt \\
& =y_{0}+hf\left( x_{1},y_{1}\right) +h\left( \frac{h}{2}\right) \left(
K\left( x_{1},y_{0},x_{0}\right) +K\left( x_{1},y_{1},x_{1}\right) \right) ,
\end{align*}

and the solution at $x_{2}$ via 
\begin{align*}
y_{2}=& \text{ }y_{1}+hf\left( x_{2},y_{2}\right)
+h\int\limits_{x_{0}}^{x_{2}}K\left( x_{2},y\left( t\right) ,t\right) dt \\
=& \text{ }y_{1}+hf\left( x_{2},y_{2}\right) \\
& \text{ }+\frac{h^{2}}{2}\left( K\left( x_{2},y_{0},x_{0}\right) +2K\left(
x_{2},y_{1},x_{1}\right) +K\left( x_{2},y_{2},x_{2}\right) \right) .
\end{align*}%
In both cases, the integral has been approximated using the composite
Trapezium Rule with values for $t$ obtained from the nodes. It must be
appreciated that, because the kernel is dependent on $x_{i}$ - the upper
limit of the integral - the integral must be calculated in its entirety at
each iteration. Hence, there will be an increasing number of terms in the
Trapezium approximation as the iteration proceeds.

In general we have%
\begin{align}
y_{i+1}=\text{ }& y_{i}+hf\left( x_{i+1},y_{i+1}\right)
+h\int\limits_{x_{0}}^{x_{i+1}}K\left( x_{i+1},y\left( t\right) ,t\right) dt
\notag \\
=\text{ }& y_{i}+hf\left( x_{i+1},y_{i+1}\right)  \label{implicit} \\
& \text{ }+\frac{h^{2}}{2}\left( \sum\limits_{j=0}^{i+1}2K\left(
x_{i+1},y_{j},x_{j}\right) -K\left( x_{i+1},y_{0},x_{0}\right) -K\left(
x_{i+1},y_{i+1},x_{i+1}\right) \right) .  \notag
\end{align}%
The presence of $y_{i+1}$ on both sides identifies this as an \textit{%
implicit} algorithm. In Appendix B we provide an insight into how this
equation can be solved for $y_{i+1}$\ using Newton's Method.

\medskip The \textit{explicit} form of (\ref{implicit}) is given by%
\begin{align}
y_{i+1}=& \text{ }y_{i}+hf\left( x_{i},y_{i}\right)  \label{explicit} \\
& \text{ }+\frac{h^{2}}{2}\left( \sum\limits_{j=0}^{i}2K\left(
x_{i},y_{j},x_{j}\right) -K\left( x_{i},y_{0},x_{0}\right) -K\left(
x_{i},y_{i},x_{i}\right) \right) .  \notag
\end{align}

\medskip

\section{\protect\medskip Stability}

The test equation we use here is%
\begin{align}
y^{\prime }\left( x\right) & =\lambda \left( y-1\right) +\gamma
\int\limits_{0}^{x}y\left( t\right) dt  \label{test eqn} \\
y\left( 0\right) & =2  \notag
\end{align}%
with solution%
\begin{align*}
y\left( x\right) & =e^{m_{1}x}+e^{m_{2}x} \\
m_{1}& =\frac{\lambda -\sqrt{\lambda ^{2}+4\gamma }}{2},\text{ \ }m_{2}=%
\frac{\lambda +\sqrt{\lambda ^{2}+4\gamma }}{2}
\end{align*}%
when $m_{1}$ and $m_{2}$ are real $\left( \lambda ^{2}+4\gamma \geqslant
0\right) $, and%
\begin{equation*}
y\left( x\right) =2e^{\frac{\lambda x}{2}}\cos \left( \frac{\sqrt{\left\vert
\lambda ^{2}+4\gamma \right\vert }}{2}x\right) 
\end{equation*}%
when $m_{1}$ and $m_{2}$ are complex $\left( \lambda ^{2}+4\gamma <0\right) .
$ This test is equation is similar to that used elsewhere \cite{brunner}\cite%
{Crisci}.

In (\ref{test eqn}), $\lambda $ is intended to represent $\partial
f/\partial y$ and $\gamma $ is intended to represent $\partial K/\partial y.$
When $\lambda >0$ and/or $\gamma >0,$ at least one of $m_{1}$ and $m_{2}$
will also be greater than zero, ensuring that the solution $y\left( x\right) 
$ is an increasing function of $x$. Note that if $\gamma >0,$ $m_{1}$ and $%
m_{2}$ cannot be complex (thus obviating the oscillatory solution), and if $%
\lambda >0$ and $m_{1}$ and $m_{2}$ are complex, then the oscillatory
solution has exponentially increasing amplitude. In all of these cases, $%
\left\vert y\left( x\right) \right\vert \rightarrow \infty $ with $x$. Only
when $\lambda \leqslant 0$ and $\gamma \leqslant 0$ do we have the case that 
$\left\vert y\left( x\right) \right\vert \rightarrow 0$ with $x$. We
acknowledge that when $\lambda =\gamma =0,$ we have the constant solution $%
y\left( x\right) =2.$

The property of \textit{stability} requires that the numerical solution to
the test equation qualitatively mimics the test solution in the sense $%
\left\vert y\left( x\right) \right\vert \rightarrow 0$ with $x.$ This means
that, at the very least, $\left\vert y_{i}\right\vert <y\left( 0\right) =2$
for all $i$, when $\lambda <0$ and $\gamma \leqslant 0,$ and $\left\vert
y_{i}\right\vert =2$ when $\lambda =\gamma =0.$

\subsection{Implicit case}

The \textit{stability function} for (\ref{implicit}) is obtained by applying
(\ref{implicit}) to the test equation:%
\begin{align*}
y_{1}& =y_{0}+h\left( \lambda \left( y_{1}-1\right) \right) +\left( \frac{%
\gamma h^{2}}{2}\right) \left( y_{0}+y_{1}\right) \\
& =2+zy_{1}-z+\frac{w}{2}\left( 2+y_{1}\right) \\
\Rightarrow y_{1}& =\frac{2-z+w}{1-z-\frac{w}{2}}=2\left( \frac{2-z+w}{2-2z-w%
}\right) \\
& \equiv P_{1}=P_{1}\left( z,w\right) ,
\end{align*}%
where we have defined $z\equiv h\lambda ,w\equiv h^{2}\gamma $ and $P_{1}$
is the stability function at $x_{1}$ - i.e. the numerical solution of the
test equation at $x_{1}$. It can be shown that%
\begin{equation*}
P_{i}=P_{i}\left( z,w\right) \equiv y_{i}=2\left( \frac{P_{i-1}-z-w+w%
\sum_{j=0}^{i-1}P_{j}}{2-2z-w}\right) ,
\end{equation*}%
where this recursion requires the definition $P_{0}\equiv y_{0}=2.$

For stability to persist in the numerical solution, we must demand%
\begin{equation}
\left\vert P_{i}\left( z,w\right) \right\vert <2  \label{P<2}
\end{equation}%
for all $i,$ whenever $\lambda <0$ and $\gamma \leqslant 0,$ and $\left\vert
y_{i}\right\vert =2$ when $\lambda =\gamma =0.$ For a given stepsize $h,$
this is implies that $\left( z,w\right) =\left( h\lambda ,h^{2}\gamma
\right) $ is located in the third quadrant in the $z-w$ system.

We define the \textit{stability region} $S$ in the third quadrant as%
\begin{equation*}
S\equiv \left\{ \left. \left( z,w\right) \right\vert \text{ }z<0,w\leqslant 0%
\text{ and }\left\vert P_{i}\left( z,w\right) \right\vert <2\text{ for all }%
i\right\} \cup \left\{ \left( 0,0\right) \right\} .
\end{equation*}%
The inclusion of $\left( z,w\right) =\left( 0,0\right) $ caters for the case 
$\lambda =\gamma =0.$

At the time of writing, we are not in possession of an analytical proof that
defines the boundary $\partial S$ of $S.$ Nevertheless, there is a useful
numerical approach that we can adopt. We consider the region 
\begin{equation*}
R_{I}\equiv \left[ -100,0\right] \times \left[ -100,0\right]
\end{equation*}%
in the $z-w$ plane. We choose a grid of $10^{6}$ points $\left( z,w\right)
\in R_{I},$ and we compute $\left\vert P_{i}\left( z,w\right) \right\vert $
for each of these points, for $i\in \left[ 1,10^{6}\right] .$ We determine,
for each $\left( z,w\right) ,$ whether or not $\left\vert P_{i}\left(
z,w\right) \right\vert >2$ for any $i\in \left[ 1,10^{6}\right] .$ If not,
we consider the point $\left( z,w\right) $ to be \textit{practically stable}%
. We have found that every single point in $R_{I}$ that we considered was
practically stable. This leads us to conjecture that all points in $R_{I}$
are stable, so that, at the very least, $S=R_{I}$. Moreover, we speculate
that $S$ is actually given by the \textit{entire} third quadrant, suggesting
unconditional stability.

Such speculation notwithstanding, the practical approach to using our
results is as follows: assume we have $\lambda \leqslant 0$ and $\gamma
\leqslant 0.$ Then the point $\left( z,w\right) =\left( \lambda ,\gamma
\right) $ is simply $\left( z,w\right) =\left( h\lambda ,h^{2}\gamma \right) 
$ with $h=1$. If such point lies within $R_{I},$ then we expect the
numerical solution, with $h=1,$ to be practically stable. It is likely to be
inaccurate, but it will be stable. Even if we reduce $h$ for the sake of
accuracy, $\left( h\lambda ,h^{2}\gamma \right) $ will still lie in $R_{I}$
and practical stability will persist. If $\left( \lambda ,\gamma \right) $
lies outside $R_{I},$ then we simply choose $h$ such that $\left( h\lambda
,h^{2}\gamma \right) \in R_{I}.$ Note that our definition of practical
stability is tied to the maximal value of $i$ that we considered $(10^{6}).$
We expect that this large enough to cater for most practical calculations;
if not, then the user can easily augment our results by checking the
condition (\ref{P<2}) for larger values for $i.$

\subsection{Explicit case}

\medskip Applying the explicit algorithm to the test equation yields%
\begin{align*}
P_{0}& =2 \\
P_{1}\left( z,w\right) & =z+2 \\
P_{2}\left( z,w\right) & =z^{2}+\frac{zw}{2}+2z+2w+2 \\
P_{i}\left( z,w\right) & \equiv 2\left( \frac{\left( 1+z+\frac{w}{2}\right)
P_{i-1}-z-w+w\sum_{j=0}^{i-2}P_{j}}{2}\right) ,\text{ \ }i\geqslant 3.
\end{align*}

Using a similar approach to that described above, we may construct a region $%
R_{E}$ of practical stability for the explicit method, and we show this
region in Figure 1. Also shown are curves, which we call $h$\textit{-paths},
for various values of $\left( \lambda ,\gamma \right) .$ The $h$-path for $%
\left( \lambda ,\gamma \right) $ is the locus of points on the diagram for $%
\left( h\lambda ,h^{2}\gamma \right) ,$ where $h\in \left[ 0,\infty \right]
. $ The $h$-path for $\left( \lambda ,\gamma \right) $ always passes through 
$\left( \lambda ,\gamma \right) $ - corresponding to $h=1$ - and always
terminates at the origin. It is clear that three of the $h$-paths shown (A,
B and C) intersect $R_{E}$ quite obviously; one of them (D), however,
intersects $R_{E}$ only near the origin. These $h$-paths indicate the effect
of reducing $h$ on the location of the point $\left( h\lambda ,h^{2}\gamma
\right) .$ Recall, for stability, it is necessary to choose $h$ such that $%
\left( h\lambda ,h^{2}\gamma \right) $ lies within the stability region, and
the $h$-path gives a clear idea of the range of values of $h$ required to
achieve that. For the $h$-path D that intersects $R_{E}$ only near the
origin, a very small value of $h$ will be needed, particularly if $\lambda $
and/or $\gamma $ are strongly negative. For the other three $h$-paths,
relatively larger values for $h$ could be tolerated to achieve stability.

\section{Error control}

We assume the algorithms have a first-order error. Consequently, we can use
Richardson extrapolation to achieve solutions of higher-order. We will not
discuss this here, and the reader is referred to our previous work for the
necessary detail \cite{prentice 2}\cite{prentice 1}. In \cite{prentice 2} we
estimate the number of nodes needed for various tolerances when using a
third-order solution. In this paper, we make similar estimates (in the next
section) with reference to a fourth-order solution. In the notation of \cite%
{prentice 2}, we use $Y_{i}^{4}-Y_{i}^{5}$ instead of $Y_{i}^{3}-Y_{i}^{5}$
for error control. Again, the reader is referred to these previous papers
for relevant detail.

\section{Examples}

\medskip We consider the same examples as used in \cite{Crisci}, modified to
suit our test equation. The parameters $\lambda $ and $\gamma $ are easily
identified. The quantities $N_{1}$ and $N_{2}$ in Table 1 refer to the
number of nodes $(N$ in (\ref{nodes})) needed to achieve tolerances of $%
\varepsilon =10^{-6}$ and $\varepsilon =10^{-12},$ respectively, using the
Richardson process described above. The stepsize can be found from $%
h=10/\left( N-1\right) .$ The implicit algorithm yielded stable solutions
for all the examples.

\begin{equation*}
\end{equation*}

\begin{center}
\renewcommand{\arraystretch}{1.7}%
\begin{tabular}{|l|c|c|cl}
\multicolumn{4}{l}{Table 1: Examples $1-3,$ implicit algorithm} &  \\ \hline
\textbf{\#} & \textbf{IDE} & $x$ & $N_{1}$ & \multicolumn{1}{|c|}{$N_{2}$}
\\ \hline
\textbf{1} & \multicolumn{1}{|l|}{$y^{\prime }=-100\left( y-1\right)
-0.1\int_{0}^{x}ydt$} & \multicolumn{1}{|l|}{$\left[ 0,10\right] $} & $1158$
& \multicolumn{1}{|c|}{$36606$} \\ \hline
\textbf{2} & \multicolumn{1}{|l|}{$y^{\prime }=-14\left( y-1\right)
-15\int_{0}^{x}ydt$} & \multicolumn{1}{|l|}{$\left[ 0,10\right] $} & $207$ & 
\multicolumn{1}{|c|}{$6519$} \\ \hline
\textbf{3} & \multicolumn{1}{|l|}{$y^{\prime }=-0.1\left( y-1\right)
-650\int_{0}^{x}ydt$} & \multicolumn{1}{|l|}{$\left[ 0,10\right] $} & $10044$
& \multicolumn{1}{|c|}{$317613$} \\ \hline
\end{tabular}%
\renewcommand{\arraystretch}{1}%
\begin{equation*}
\end{equation*}
\end{center}

Results for the explicit algorithm are shown in Table 2. Here, $N$ indicates
the least number of nodes required for a stable solution, $h_{s}$ is the
corresponding stepsize and $\left( z,w\right) =\left( h_{s}\lambda
,h_{s}^{2}\gamma \right) .$

\begin{center}
\renewcommand{\arraystretch}{1.7}%
\begin{tabular}{|l|c|c|cc}
\multicolumn{4}{l}{Table 2: explicit algorithm} &  \\ \cline{1-4}
\textbf{\#} & $N$ & \multicolumn{1}{|c|}{$h_{s}$} & $\left( z,w\right) $ & 
\multicolumn{1}{|c}{} \\ \cline{1-4}
\textbf{1} & \multicolumn{1}{|l|}{$505$} & \multicolumn{1}{|l|}{$0.0198$} & $%
\left( -1.98,4\times 10^{-5}\right) $ & \multicolumn{1}{|c}{} \\ \cline{1-4}
\textbf{2} & \multicolumn{1}{|l|}{$72$} & \multicolumn{1}{|l|}{$0.1408$} & $%
\left( 1.97,0.3\right) $ & \multicolumn{1}{|c}{} \\ \cline{1-4}
\textbf{3} & \multicolumn{1}{|l|}{$9501$} & \multicolumn{1}{|l|}{$0.0011$} & 
$\left( -1.1\times 10^{-5},7.8\times 10^{-4}\right) $ & \multicolumn{1}{|c}{}
\\ \cline{1-4}
\end{tabular}%
\renewcommand{\arraystretch}{1}
\end{center}

Note the very small stepsize needed for \#3. This case corresponds to the $h$%
-path D in Figure 1 that is close to the vertical axis, and intersects the
stability region only near the origin.

\section{Conclusion}

We have extended earlier work on the numerical solution of Volterra
integro-differential equations by providing implicit and explicit versions
of an Euler-type algorithm for a kernel of general form. We have conducted a
stability analysis. We have found that the implicit method is stable over a
large region in the third quadrant of the stability space (and very possibly
over the entire third quadrant), while the explicit method is stable over a
relatively small region near the origin. Numerical examples with a suitable
test equation support these findings.

\section{Appendix A}

\subsection{Systems}

When $n=2$ in (\ref{problem}), we have the system 
\begin{equation*}
\left[ 
\begin{array}{c}
y_{1}^{\prime } \\ 
y_{2}^{\prime }%
\end{array}%
\right] =\left[ 
\begin{array}{l}
y_{2} \\ 
f\left( x,y_{1}\right) +\int\limits_{x_{0}}^{x}K\left( x,y_{1}\left(
t\right) ,t\right) dt%
\end{array}%
\right]
\end{equation*}%
and when $n=3,$ we have%
\begin{equation*}
\left[ 
\begin{array}{c}
y_{1}^{\prime } \\ 
y_{2}^{\prime } \\ 
y_{3}^{\prime }%
\end{array}%
\right] =\left[ 
\begin{array}{l}
y_{2} \\ 
y_{3} \\ 
f\left( x,y_{1}\right) +\int\limits_{x_{0}}^{x}K\left( x,y_{1}\left(
t\right) ,t\right) dt%
\end{array}%
\right] ,
\end{equation*}%
The initial values $y_{1}\left( x_{0}\right) $ and $y_{2}\left( x_{0}\right)
=y_{1}^{\prime }\left( x_{0}\right) $ must be specified for the first
system, and $y_{1}\left( x_{0}\right) ,y_{2}\left( x_{0}\right)
=y_{1}^{\prime }\left( x_{0}\right) $ and $y_{3}\left( x_{0}\right)
=y_{2}^{\prime \prime }\left( x_{0}\right) $ must be specified for the
second system. Of course, analogous expressions exist for $n>3$.

The general first-order system has the form

\begin{align*}
\left[ 
\begin{array}{c}
y_{1}^{\prime } \\ 
\vdots \\ 
y_{m}^{\prime }%
\end{array}%
\right] & =\left[ 
\begin{array}{l}
f_{1}\left( x,y_{1},\ldots ,y_{m}\right) +\int\limits_{x_{0}}^{x}K_{1}\left(
x,y_{1},\ldots ,y_{m}\right) dt \\ 
\multicolumn{1}{c}{\vdots} \\ 
f_{m}\left( x,y_{1},\ldots ,y_{m}\right) +\int\limits_{x_{0}}^{x}K_{m}\left(
x,y_{1},\ldots ,y_{m}\right) dt%
\end{array}%
\right] \\
\Rightarrow \left[ 
\begin{array}{c}
y_{1,i+1} \\ 
\vdots \\ 
y_{m,i+1}%
\end{array}%
\right] & =\left[ 
\begin{array}{l}
y_{1,i}+hG_{1} \\ 
\multicolumn{1}{c}{\vdots} \\ 
y_{m,i}+hG_{m}%
\end{array}%
\right] ,
\end{align*}%
where the forms of $G_{1}$ and $G_{m}$ can be inferred. For this system, we
can find the Jacobians%
\begin{equation*}
J_{f}\equiv \left[ 
\begin{array}{ccc}
\frac{\partial f_{1}}{\partial y_{1}} & \cdots & \frac{\partial f_{1}}{%
\partial y_{m}} \\ 
\vdots & \ddots & \vdots \\ 
\frac{\partial f_{m}}{\partial y_{1}} & \cdots & \frac{\partial f_{m}}{%
\partial y_{m}}%
\end{array}%
\right] \text{ and }J_{K}\equiv \left[ 
\begin{array}{ccc}
\frac{\partial K_{1}}{\partial y_{1}} & \cdots & \frac{\partial fK_{1}}{%
\partial y_{m}} \\ 
\vdots & \ddots & \vdots \\ 
\frac{\partial K_{m}}{\partial y_{1}} & \cdots & \frac{\partial K_{m}}{%
\partial y_{m}}%
\end{array}%
\right] .
\end{equation*}%
Let $\left\{ \lambda _{1},\ldots ,\lambda _{m}\right\} $ denote the
eigenvalues of $J_{f},$ and $\left\{ \gamma _{1},\ldots ,\gamma _{m}\right\} 
$ denote the eigenvalues of $J_{K}.$ We form the pairs $\left( \lambda
_{j},\gamma _{j}\right) .$ If $\lambda _{j}\leqslant 0$ and $\gamma
_{j}\leqslant 0$, or $\func{Re}\lambda _{j}\leqslant 0$ and $\func{Re}\gamma
_{j}\leqslant 0$ if the eigenvalues are complex, we write%
\begin{equation*}
\left( z,w\right) =\left( h\lambda _{j},h^{2}\gamma _{j}\right) =\left( -h%
\sqrt{\lambda _{j}^{\ast }\lambda _{j}},-h^{2}\sqrt{\gamma _{j}^{\ast
}\gamma _{j}}\right) .
\end{equation*}

\section{Appendix B}

To solve (\ref{implicit}) for $y_{i+1}$ we employ Newton's Method in the form%
\begin{equation*}
y_{i+1}^{k+1}=y_{i+1}^{k}-\frac{F\left( y_{i+1}^{k}\right) }{F^{\prime
}\left( y_{i+1}^{k}\right) },
\end{equation*}%
where 
\begin{align*}
F\left( y_{i+1}\right) \equiv \text{ }& y_{i+1}-y_{i}-hf\left(
x_{i+1},y_{i+1}\right) \\
& \text{ }-\frac{h^{2}}{2}\left( \sum\limits_{j=0}^{i+1}2K\left(
x_{i+1},y_{j},x_{j}\right) -K\left( x_{i+1},y_{0},x_{0}\right) -K\left(
x_{i+1},y_{i+1},x_{i+1}\right) \right)
\end{align*}%
and%
\begin{equation*}
F^{\prime }\left( y_{i+1}\right) =\frac{dF\left( y_{i+1}\right) }{dy_{i+1}}%
=1-h\frac{df\left( x_{i+1},y_{i+1}\right) }{dy_{i+1}}-\frac{h^{2}}{2}\frac{%
dK\left( x_{i+1},y_{i+1},x_{i+1}\right) }{dy_{i+1}},
\end{equation*}%
with initial value $y_{i+1}^{0}=y_{i}.$ We would generally expect
second-order convergence for this method.

If we use%
\begin{equation*}
y_{i+1}^{k+1}=y_{i+1}^{k}-\frac{F\left( y_{i+1}^{k}\right) }{F^{\prime
}\left( y_{i+1}^{k}\right) }-\frac{F^{2}\left( y_{i+1}^{k}\right) F^{\prime
\prime }\left( y_{i+1}^{k}\right) }{2\left( F^{\prime }\left(
y_{i+1}^{k}\right) \right) ^{3}}
\end{equation*}%
where%
\begin{equation*}
F^{\prime \prime }\left( y_{i+1}\right) =\frac{d}{dy_{i+1}}\frac{dF\left(
y_{i+1}\right) }{dy_{i+1}},
\end{equation*}%
we may expect third-order convergence, which may improve the efficiency of
the algorithm.

\end{document}